\documentstyle[12pt]{article}
\begin{document}

\author{S.V. Ludkowski}

\title{Separability of metagroup algebras.}

\date{21 September 2018}
\maketitle

\begin{abstract}
For a class of nonassociative metagroup algebras their separability
is investigated. For this purpose the cohomology theory on them is
utilized. Conditions are found under which nonassociative metagroup
algebras are separable. Algebras satisfying these conditions are
described. \footnote{key words and phrases: algebra; nonassociative;
separable; ideal;
cohomology 
address: Dep. Appl. Mathematics, Moscow State Techn. Univ. MIREA,
av. Vernadsky 78, Moscow, 119454, Russia, sludkowski@mail.ru }

\end{abstract}

\section{Introduction.}
Associative separable algebras play an important role and have found
many-sided application (see, for example,
\cite{bourbalgch123b,flosaam11,gesagmj78,mapesaam09,montsmith,
saoystha,pierceb,rumsaal94} and references therein). Studies of
their structure are based on cohomology theory. On the other hand,
cohomology theory of associative algebras was investigated by
Hochschild and other authors
\cite{cartaneilenbb56,hochschild46,pommb}, but it is not applicable
to nonassociative algebras. Cohomology theory of group algebras is
an important and great part of algebraic topology. It is worth to
mention that nonassociative algebras with some identities in them
found many-sided applications in physics, noncommutative geometry,
quantum field theory, PDEs and other sciences (see
\cite{dickson,girardgb,guertzeb,kansol,krausshb,ludwrgrijmgta}-
\cite{lustsdoadaca11,schaeferb,serodaaca07} and references therein).
\par An extensive area of investigations of
PDE intersects with cohomologies and deformed cohomologies
\cite{pommb}. Therefore, it is important to develop this area over
octonions, Cayley-Dickson algebras and more general metagroup
algebras. Some results in this area are presented in
\cite{ludcohmalosal17}.
\par This article is devoted to a separability of nonassociative
metagroup algebras. Conditions are found under which they are
separable. Algebras satisfying these conditions are described.
\par A formula "$(m)$" within the same subsection "n" and section $"k"$
is referred as "$(m)$"; in another subsection within the same
section as "$n(m)$"; in another section as "$k.n(m)$". A theorem,
lemma, etc. number "n" within the same section "k" is referred as
"n", in another section as "k.n".

\par All main results of this paper are obtained for the first time.
They can be used for further studies of nonassociative algebras,
their cohomologies, algebraic geometry, PDEs, their applications in
the sciences, etc.

\section{Separable nonassociative algebras.}
To avoid misunderstandings we first give our definitions and
notations.
\par {\bf 1. Definition.} Let \par $(1)$ $A$ be a nonassociative metagroup algebra over a commutative
associative unital ring $\cal T$ with a metagroup $G$ such that
$G\cap {\cal T} = {\bf \Psi }$. \par A $G$-graded $A$-module $P$
(see subsection 2.7 in \cite{ludcohmalosal17} also) is called
projective if it is isomorphic with a direct additive of a free
$G$-graded $A$-module. The metagroup algebra $A$ is called separable
if it is a projective $G$-graded $A^e$-module.
\par One puts $\mu (z)=1_Az$ for each $z\in A^e$, where $A$ is
considered as the $G$-graded right $A^e$-module.

\par {\bf 2. Proposition.} {\it Suppose that $A$ is a nonassociative algebra
satisfying condition 1$(1)$. Then the following conditions are
equivalent:
\par $(1)$ $A$ is separable;
\par $(2)$ the exact sequence $0\rightarrow Ker ~ \mu \rightarrow A^e\mbox{ }_{\overrightarrow{\mu }} A\rightarrow 0$
splits;
\par $(3)$ an element $b\in A^e$ exists such that $\mu (b)=1_A$ and
$xb=bx$ and $b(xy)=(bx)y$ and $(xb)y=x(by)$ and $(xy)b=x(yb)$ for
all $x$ and $y$ in $A$, where $A^e$ is considered as the $G$-graded
two-sided $A$-module.}
\par {\bf Proof.} The  implication $(1)\Rightarrow (2)$ is evident.
\par $(2)\Rightarrow (3)$. If the exact sequence $(2)$ splits, then $A^e$ as
the $A^e$-module is isomorphic with $A\oplus ker (\mu )$. Therefore,
$A$ is separable. The sequence $(2)$ splits if and only if there
exists $p\in Hom_{A^e}(A,A^e)$ such that $\mu p=id_A$. With this
homomorphism $p$ put $b=p(1_A)$. Then
$(xb)y=(xp(1_A))y=p(x1_A)y=p(x(1_Ay))=p((xy)1_A)=(xy)b$, hence $\mu
(b)=\mu p(1_A)=1_A$ and $xb=xp(1_A)=p(x1_A)=p(1_Ax)=p(1_A)x=bx$.
Thus $(3)$ is valid.
\par $(3)\Rightarrow (1)$. Suppose that condition $(3)$ is
fulfilled. Then a mapping $p : A\to A^e$ exists such that $p(x)=bx$.
The element $b$ has the decomposition $b=\sum_j b_j g_j$ with
$g_j=g_{j,1}\otimes g_{j,2}$, where $g_{j,1}\in G$ and $g_{j,2}\in
G^{op}$ and $b_j\in \cal T$ for each $j$.  Therefore, using
conditions $(3)$ and 2.4$(1-4)$ \cite{ludcohmalosal17} we infer that
\par $p(xy)=\sum_j \sum_k \sum_l b_jg_j((c_kx_k)(d_ly_l)) = \sum_j
b_j(g_jx)y= (bx)y= p(x)y$ and
\par $p(yx)=(by)x=(yb)x=y(bx)=yp(x)$ for each $x$ and $y\in A$, where
$x=\sum_k c_kx_k$ and $y=\sum_ld_ly_l$ with $x_k$ and $y_l$ in $G$,
$c_k$ and $d_l$ in ${\cal T}$ for each $k$ and $l$. Thus $p\in
Hom_{A^e}(A,A^e)$. Moreover, $\mu (p(x))=\mu (bx)=\mu (b)x=1_Ax=x$
for each $x\in A$, consequently, the exact sequence $(2)$ splits.

\par {\bf 3. Definition.} An element $b\in A^e$ fulfilling condition
$(3)$ is called a separating idempotent of an algebra $A$.

\par {\bf 4. Lemma.} {\it Let $A$ be a nonassociative algebra
satisfying condition 1$(1)$. Let also $M$ be a two-sided $A$-module.
\par $(1)$. If $p\in Hom_{A^e}(ker (\mu ), M)$ and $\kappa : A\to A^e$ with
$\kappa (x)=x\otimes 1-1\otimes x$ for each $x\in A$, then $p\kappa
$ is a derivation of $A$ with values in $M$. \par $(2)$. A mapping
$\chi : p\mapsto  p\kappa $ is an isomorphism of $Hom_{A^e}(ker (\mu
), M)$ onto $Z^1_{\cal T}(A,M)$.
\par $(3)$. $\chi ^{-1}(B^1_{\cal T}(A,M))= \{ \psi |_{\ker (\mu )}:
\psi \in Hom_{A^e}(A^e,M) \} $.}
\par {\bf Proof.} $(1)$. Since $\mu \kappa =0$, then $Im (\kappa )\subseteq
ker (\mu )$. By virtue of Proposition 2.12 \cite{ludcohmalosal17}
$\mu \kappa $ is the derivation having also properties $(2)$ and
$(3)$.

\par {\bf 5. Proposition.} {\it Suppose that $A$ is a nontrivial
nonassociative algebra satisfying condition 1$(1)$. Then $H^1_{\cal
T}(A,M)=0$ for each two-sided $A$-module $M$ if and only if $A$ is a
separable ${\cal T}$-algebra.}
\par {\bf Proof.}  In view of Proposition 2 the algebra $A$ is
separable if and only if the exact sequence 2$(2)$ splits. That is a
homomorphism $h$ exists $h\in Hom _{A^e}(A, ker (\mu ))$ such that
its restriction $h|_{ker (\mu )}$ is the identity mapping.
Therefore, if $H^1_{\cal T}(A, ker (\mu ))=0$, then the algebra $A$
is separable due to Lemma 4. \par Vice versa if a homomorphism $h\in
Hom _{A^e}(A^e, ker (\mu ))$ exists with $h|_{ker (\mu )}=id$, then
each $p\in Hom _{A^e}(ker (\mu ), M)$ has the form $f|_{ker (\mu )}$
with $f=ph\in Hom _{A^e}(A^e, M)$. By virtue of Lemma 4 $Z^1_{\cal
T}(A,M)=B^1_{\cal T}(A,M)$ for each two-sided $A$-module $M$.

\par {\bf 6. Theorem.} {\it Let a noncommutative algebra $A$ fulfill
condition 1$(1)$ and
\par $(1)$ $Dim (A/J(A))\le 1$ and
\par $(2)$ $A/J(A)$ is projective as the $\cal T$-module and
\par $(3)$ $J(A)^k=0$ for some $k\ge 1$, where $J(A)$ denotes the radical of $A$.
\par Then a subalgebra $D$ in $A$ exists such that $A=D\oplus J(A)$
as $\cal T$-modules and $A/J(A)$ is isomorphic with $D$ as the
algebra.}
\par {\bf Proof.} For $k=1$ we get $A=D$. \par For $k=2$ a natural
projection $\pi : A\to A/J$ exists, where $J=J(A)=rad (A_A)$, since
$J^2=0$. The algebra $A$ is $G$-graded and ${\cal T}\subset Z(A)$,
hence $rad ((A_e)_{A_e}) \subseteq (rad (A_A))_e$, where $e$ is the
unit element of $G$. In view of conditions 2.7$(1-3)$ in
\cite{ludcohmalosal17} $J$ is the two-sided ideal in $A$ and
$J_r^m=J_l^m$ for each positive integer $m$, where $J_l^1=J$,
$J_r^1=J$, $J_l^{m+1}=JJ_l^m$ and $J_r^{m+1}=J_r^mJ$. Conditions
2.1$(9)$ and 2.7$(1-3)$ in \cite{ludcohmalosal17} imply that $A/J$
is also $G$-graded, since ${\cal T}\subset Z(A)$. \par By condition
$(2)$ the $\cal T$-module $A/J$ is projective, consequently, an
exact splitting sequence of $\cal T$-modules exists $0\to J\to A\to
A/J\to 0$. Thus a homomorphism $\kappa : A/J\to A$ of $\cal
T$-modules exists such that $\pi \kappa =id$ on $A/J$. For any two
elements $x$ and $y$ in $A/J$ put $\Phi (x,y)=\kappa (xy)-\kappa
(x)\kappa (y)$. Therefore, $\pi \Phi (x,y)=\pi \kappa (xy)-\pi
(\kappa (x)\kappa (y))=xy-xy=0$, since $\pi $ is the algebra
homomorphism and $\pi \kappa =id$. Thus $\Phi (x,y)\in ker (\pi
)=J$.
 One has by the definition that \par $Dim (A/J)=\sup \{ n: \exists  \mbox{  two-sided }
A/J \mbox{-module } M ~ ~ H^n_{\cal T}(A/J,M)\ne 0 \} $.
\par Then put $ux=u\kappa (x)$ and $xu=\kappa (x)u$ to be the right
and left actions of $A/J$ on $J$. Since $\kappa $ is the
homomorphism of $\cal T$-modules and ${\cal T}\subseteq Z(A)$, then
for each pure states $x$, $y$ and $u$: $(xy)u-{\sf t}_3x(yu)=\kappa
(xy)u-(\kappa (x)\kappa (y))u=\Phi (x,y)u\in J^2=0$, where ${\sf
t}_3={\sf t}_3(x,y,u)$. Then $u(xy)-{\sf t}_3^{-1}(ux)y=u\kappa
(xy)-u(\kappa (x)\kappa (y)) = u\Phi (x,y)\in J^2=0$, where ${\sf
t}_3={\sf t}_3(u,x,y)$. Thus $J$ has the structure of the two-sided
$A/J$-module.
\par Evidently, $\Phi $ is $\cal T$-bilinear. Then for every
pure states $x$, $y$ and $z$ in $A/J$: \par $(\delta ^2\Phi
)(x,y,z)={\sf t}_3x(\kappa (yz)-\kappa (y)\kappa (z))-(\kappa
((xy)z)-\kappa (xy)\kappa (z))+{\sf t}_3(\kappa (x(yz))-\kappa
(x)\kappa (yz))-(\kappa (xy)-\kappa (x)\kappa (y))z$\par $={\sf
t}_3\kappa (x)\kappa (yz)- {\sf t}_3\kappa (x)(\kappa (y)\kappa
(z))- \kappa ((xy)z)+\kappa (xy)\kappa (z) +{\sf t}_3\kappa
(x(yz))-{\sf t}_3\kappa (x)\kappa (yz)-\kappa (xy)\kappa (z)+(\kappa
(x)\kappa (y))\kappa (z)=0$, \\ consequently, $\Phi \in B^2_{\cal
T}(A/J,J)$, where ${\sf t}_3={\sf t}_3(x,y,z)$. Thus by the $\cal
T$-linearity a homomorphism $h$ in $Hom_{\cal T}(A/J,J)$ exists
possessing the property $\Phi (x,y)=xh(y)-h(xy)-h(x)y$ for each $x$
and $y$ in $A/J$.
\par Let now $p=\kappa +h\in Hom_{\cal T}(A/J,J)$, consequently,
$\pi p=\pi \kappa =id|_{A/J}$, since $\pi (J)=0$. This implies that
$p(xy)-p(x)p(y)=0$ for each $x$ and $y$ in $A/J$, since $\kappa
(xy)-\kappa (x)\kappa (y)=\Phi (x,y)=xh(y)-h(xy)+h(x)y$ and
$h(x)h(y)\in J^2=0$. Since $p(1_{A/J})-1_A\in J$, then
$(p(1)-1)^2=1-p(1)$. Therefore, $p$ is the algebra homomorphism.
This implies that $D=Im (p)$ is the subalgebra in $A$ such that
$A=D\oplus J$.
\par Let now $k>2$ and this theorem is proven for $1,...,k-1$.
Put $A_1=A/J^2$, then $J/J^2$ is the  two-sided ideal in $A_1$ and
$A_1/(J/J^2)$ is isomorphic with $A/J$, also $(J/J^2)^2=0$. Thus
$J(A_1)=J/J^2$ and $A_1$ satisfies conditions $(1-3)$ of this
theorem and is $G$-graded, since $A$ and $J$ are $G$-graded and
${\cal T}\subset Z(G)$ due to conditions 2.1$(9)$ and 2.7$(1-3)$ in
\cite{ludcohmalosal17}.
\par From the proof for $k=2$ we get that a subalgebra $D_1$ in $A_1$
exists such that $A_1=D_1\oplus J/J^2$. Consider a subalgebra $E$ in
$D$ such that $E\cap J=J^2$ and $D_1=E/J^2$. Then $E/J$ is
isomorphic with $E/(E\cap J) \approx (E+J)/J=A/J$. Moreover,
$(J^2)^{k-1}=J^{k+k-2}\subseteq J^k=0$, hence $J(E)=J^2$. Thus the
algebra $E$ fulfills conditions $(1-3)$ of this theorem, is
$G$-graded and $J(E)^{k-1}=0$. \par By the induction supposition a
subalgebra $F$ in $E$ exists such that $E=F\oplus J^2$,
consequently, $F+J=E+J=A$ and $F\cap J=F\cap E\cap J=F\cap J^2=0$.
Thus $A=F\oplus J$.
\par {\bf 7. Proposition.} {\it Suppose that conditions of Theorem 6
are satisfied and $(1)$ takes the form $Dim (A/J(A))=0$. Then for
any two $G$-graded subalgebras $B$ and $C$ in $A$ such that
$A=B\oplus J(A)$ and $A=C\oplus J(A)$ an element $v\in J(A)$ exists
for which $(1-v)C=B(1-v)$ such that $(1-v)$ has a right inverse and
a left inverse.}
\par {\bf Proof.} Let $q: A\to B$ and $r: A\to C$ be the canonical projections
induced by the decompositions $A=B\oplus J$ and $A=C\oplus J$, where
$J=J(A)$. Then $p \pi =q$ and $s \pi =r$, where $\pi : A\to A/J$ is
the quotient homomorphism, $p: A/J\to A$ and $s: A/J\to C$ are
homomorphisms as in the proof of Theorem 6, since $q$ and $r$ are
homomorphisms of algebras. Put $w(x)=p(x)-s(x)$ for each $x\in A/J$,
$w: A/J\to J$. Then $\pi (w\pi )=\pi (p\pi )-\pi (s\pi )=\pi q-\pi
r=\pi (id_A-r)-\pi (id_A-q)=0$, since $Im (id_A-q)=Im (id_A-r)=J=ker
(\pi )$. Therefore, $Im (w)=Im (w\pi )\subseteq J$, hence $w\in
Hom_{\cal T} (A/J,J)$. Then
$w(xy)=p(xy)-s(xy)=p(x)(p(y)-s(y))+(p(x)-s(x))s(y)=xw(y)+w(x)y$,
consequently, $w$ is the derivation of the algebra $A/J$ with values
in the two-sided $A$-module $A/J$ (see also the proof of Theorem 6).
Since $Dim (A/J)=0$, then $w$ is the inner derivation by Proposition
2.12 in \cite{ludcohmalosal17}. Thus an element $v\in J$ exists for
which $w(x)=xv-vx$ for each $x\in A/J$. This implies that
$p(x)(1-v)=(1-v)s(x)$ for each $x\in A/J$. The element $(1-v)$ has a
right inverse and a left inverse, since $J^k=0$ implies $v_l^k=0$
and $v_r^k=0$, where $v_l^1=v$, $v_r^1=v$, $v_l^{m+1}=vv_l^m$ and
$v_r^{m+1}=v_r^mv$ for each positive integer $m$. Therefore,
$B(1-v)=p(A/J)(1-v)=(1-v)s(A/J)=(1-v)C$.

\end{document}